\documentclass[preprint,12pt]{elsarticle}
\journal{J.}% Mathematical Analysis and Applications}
\usepackage{dsfont}
\usepackage{amsfonts}
\usepackage{amsmath}
\usepackage{epsfig}
\usepackage{amssymb}
\begin{document}
\def\R{\mathbb{R}}
\def\C{\mathbb{C}}
\def\Z{\mathbb{Z}}
\def\N{\mathbb{N}}
\def\Q{\mathbb{Q}}
\def\D{\mathbb{D}}
\def\T{\mathbb{T}}
\def\hb{\hfil \break}
\def\ni{\noindent}
\def\i{\indent}
\def\a{\alpha}
\def\b{\beta}
\def\e{\epsilon}
\def\d{\delta}
\def\De{\Delta}
\def\g{\gamma}
\def\qq{\qquad}
\def\L{\Lambda}
\def\E{\cal E}
\def\G{\Gamma}
\def\F{\cal F}
\def\K{\cal K}
\def\A{\cal A}
\def\B{\cal B}
\def\M{\cal M}
\def\P{\cal P}
\def\Om{\Omega}
\def\om{\omega}
\def\s{\sigma}
\def\t{\theta}
\def\th{\theta}
\def\Th{\Theta}
\def\z{\zeta}
\def\p{\phi}
\def\m{\mu}
\def\n{\nu}
\def\l{\lambda}
\def\Si{\Sigma}
\def\q{\quad}
\def\qq{\qquad}
\def\half{\frac{1}{2}}
\def\hb{\hfil \break}
\def\half{\frac{1}{2}}
\def\pa{\partial}
\def\r{\rho}
\bibliographystyle{elsarticle-harv}
\begin{frontmatter}
\title{Logarithmic moving averages}
\author{N. H. Bingham\footnote{Department of Mathematics, Imperial College London, 180 Queens Gate, London, SW7 2BZ, UK; Email: N. Bingham@imperial.ac.uk} and Bujar Gashi\footnote{Department of Mathematical Sciences, The University of Liverpool, Liverpool, L69 7ZL, UK; Email: Bujar.Gashi@liverpool.ac.uk}}

\begin{abstract}
 We introduce a moving average summability method, which is proved to be equivalent with the logarithmic $\ell$-method. Several equivalence and Tauberian theorems are given. A strong law of large numbers is also proved.
\end{abstract}


\end{frontmatter}

\numberwithin{equation}{section}
\newtheorem{proof}{Proof}
\newtheorem{definition}{Definition}
\newtheorem{theorem}{Theorem}
\newtheorem{lemma}{Lemma}
\newtheorem{proposition}{Proposition}
\newtheorem{remark}{Remark}
\newtheorem{corollary}{Corollary}
\newtheorem{assumption}{Asssumption}

\section{Introduction} The logarithmic methods of summation $\ell$ and $L$ are classical (see, for example, Hardy ~\cite{Har}, Ishiguro ~\cite{Ish0} - \cite{Ish3}). Let $\{s_n\}_{n=0}^\infty$ be a sequence of real numbers. The sequence is summable to $s$ by the logarithmic $\ell$-method, written $s_n\rightarrow s$ $(\ell)$, if
\begin{eqnarray}
t_n := \frac{1}{\log n}\sum_{i=0}^n\frac{s_i}{i+1}\rightarrow s, \quad (n\rightarrow\infty)\label{l method}
\end{eqnarray}
(we write $\ell_x$ when the limit is taken through a continuous variable).\\

The sequence is summable to $s$ by the logarithmic $L$-method, written $s_n\rightarrow s$ $(L)$, if
\begin{eqnarray}
\frac{1}{- \log(1-x)}\sum_{i=0}^\infty\frac{s_i}{i+1}x^{i+1}\rightarrow s\quad (x\uparrow 1).\label{L method}
\end{eqnarray}

Here we introduce a certain delayed (deferred) summability method. For $\lambda > 1$, the sequence $\{s_n\}_{n=0}^\infty$ is summable by the {\it logarithmic moving average}, $s_n\rightarrow s$ ($\mathcal{L}(\lambda)$), if
\begin{eqnarray}
\frac{1}{\log n}\sum_{n^{1/\lambda}<i\leq n}\frac{s_i}{i+1}\rightarrow (1-\lambda^{-1})s,\quad (n\rightarrow\infty)\label{Ln}
\end{eqnarray}
(we write ${\cal L}_x(\lambda)$ if the limit is taken through a continuous variable $x$.)
%\begin{eqnarray}
%\frac{1}{\log x}\sum_{x^{1/\lambda}<k\leq x}\frac{s_k}{k+1}\rightarrow (1-\lambda^{-1})s,\quad (x\rightarrow\infty),\label{Lx}
%\end{eqnarray}
%we use the notation $s_n\rightarrow s$ ($\mathcal{L}_x(\lambda)$).
\section{Equivalence and Tauberian theorems}
\subsection{Results}
\begin{theorem}\label{e1}$\ell\Leftrightarrow \mathcal{L}(\lambda)$ for some (all) $\lambda>1$.
\end{theorem}
With $\{a_n\}_{n=0}^\infty$ defined by $s_n=\sum_{k=0}^na_n$, the Riesz (typical) mean $R(\log n)$ (of order 1) is defined as
\begin{eqnarray}
\frac{1}{x}\int_0^x\{\sum_{k:\log(k+1)<y}a_k\} dy,\quad (x\rightarrow\infty).\nonumber
\end{eqnarray}
In view of $R(\log n) \Leftrightarrow \ell$ (Hardy ~\cite{Har} Th. 37; see also \S 5.16), this gives
\begin{corollary}\label{e3} $R(\log n)\Leftrightarrow\mathcal{L}(\lambda)$ for some (all) $\lambda>1$.
\end{corollary}
Note that $R(\log n)$ involves a {\it continuous} limit, but $\mathcal{L}(\lambda)$ a {\it discrete} one.  This equivalence between discrete and continuous limits is a consequence of {\it uniformity}, as in Theorem 2 below.
\begin{theorem}\label{u} If (\ref{Ln}) holds for all $\lambda>1$, then it holds uniformly on compact $\lambda$-sets in $(1,\infty)$.
%(ii) The convergence in (\ref{Ln}) holds for all $\lambda>1$ if and only if the convergence in (\ref{Lx}) also holds for all $\lambda>1$.
\end{theorem}
\begin{corollary}\label{e4} $s_n\rightarrow s$ $(\ell_x)$ if and only if $s_n\rightarrow s$ ($\mathcal{L}_x(\lambda)$) for all $\lambda>1$.
\end{corollary}
\begin{theorem}\label{e6} Let $U(x):=\sum_{0\leq i\leq x}s_i(i+1)^{-1}$. The following statements are equivalent:\\

(i) $U(x)=U_1(x)-U_2(x)$, with $U_2(x)$ non-decreasing and $U_1(x)$ satisfying
\begin{eqnarray}
\lim_{x\rightarrow\infty}[U_1(x)-U_1(x^{1/\lambda})](\log x)^{-1}=s(1-\lambda^{-1}),\quad \forall \lambda>1,\nonumber
\end{eqnarray}
(ii) $\liminf_{\alpha\downarrow1} \limsup_{x\rightarrow\infty}\sup_{\theta\in[1,\alpha]}[U(x)-U(x^{1/\theta})](\log x)^{-1}<\infty$.
\end{theorem}

\begin{corollary}\label{e7} If $s_n\rightarrow s$ ($\ell$), then statement (ii) of Theorem \ref{e6} holds.
\end{corollary}

The Abelian result that $\ell\Rightarrow L$ is proved in~\cite{Ish1}. The simplest Tauberian condition for $L\Rightarrow \ell$, and thus $\ell\Leftrightarrow L$, is $s_n=O_L(1)$ as proved in~\cite{Ish2}. We next give a Tauberian theorem establishing the equivalence between $\ell$ and $L$ methods under a one-sided Tauberian condition of best possible character.
\begin{theorem}\label{t1} We have $\ell\Leftrightarrow L$ if and only if
\begin{eqnarray}
\lim_{\lambda\downarrow 1}\liminf_{n\rightarrow \infty}\min_{n\leq m\leq \lambda n}\frac{1}{\log n}\sum_{n\leq i\leq m}\frac{s_i}{i+1}\geq 0.\label{one sided}
\end{eqnarray}
\end{theorem}

Let $P$ be a probability law on the non-negative integers with mean $\mu>0$, variance $\sigma^2$, and finite third moment. Consider the sequence of independent random variables $Y,Y_0,Y_1,...,$ with law $P$, and let $S_0=0$, $S_n=\sum_{k=0}^nY_k$. We write $s_n\rightarrow s$ $(P)$ if
\begin{eqnarray}
\sum_{i=0}^\infty s_iP(S_n=i)\rightarrow s,\quad (n\rightarrow \infty).\nonumber
\end{eqnarray}
This is the $P$-{\it method} or {\it random-walk method} of summation introduced in~\cite{Bin1}, ~\cite{Bin3},~\cite{Bin4}.

% The following two Tauberian  theorems establish relations between the methods $P$ and $\ell$ (equivalently, ${\cal L}(\lambda)$).
The Borel method and its relatives (Euler, Valiron and circle
methods, etc.) were studied by Hardy, Littlewood and others for their
applications to power series, as they sum more power series than the
Ces\`aro methods.  In this regard, Hardy and Littlewood obtained a
number of Tauberian theorems from Borel-type to Ces\`aro-type methods.
These were extended to random-walk ($P$-) methods by the first author.
In the next result, such results are extended to the logarithmic method
(as Ces\`aro convergence implies logarithmic convergence; see e.g.
Ishiguro~\cite{Ish1}); so too is the first author's result in the converse
direction, again using a result of Ishiguro~\cite{Ish1}.
%and the $l$-method (or equivalently the $\mathcal{L}(\lambda)$-method)
\begin{theorem}\label{t2} (i) If $s_n=O(\sqrt{n})$, then $s_n\rightarrow s$ $(P)$ implies $s_n\rightarrow s$ $(\ell)$.\\

(ii) Let $P$ have a finite $k$-th moment $(k=3,4,...)$, and $0\leq r\leq (k-2)/2$. $s_n\rightarrow s$ $(P)$ implies $s_n\rightarrow s$ $(L)$ if
\begin{eqnarray}
\lim_{\delta\downarrow0} \liminf_{n\rightarrow\infty}\min_{n\leq m<n+\delta\sqrt{n}}(s_m-s_n)n^{-r}\geq 0.\nonumber
\end{eqnarray}
If in addition $(\ref{one sided})$ holds or $k=3$, then $s_n\rightarrow s$ $(\ell)$.\\

(iii) If there exist $\varepsilon_n\rightarrow 0$ such that
\begin{eqnarray}
\frac{1}{\log n}\sum_{i=0}^n\frac{s_i+\varepsilon_i}{i+1}=s+o(1/n^{1/2}\log n),\label{P}
\end{eqnarray}
then $s_n\rightarrow s$ $(P)$.
%$$
%then $s_n$ is summable by the $P$-method.
%\begin{eqnarray}
%t_n\equiv\frac{1}{\log n}\sum_{i=0}^n\frac{s_n}{i+1}.\nonumber
%\end{eqnarray}
%$t_n= s+o(1/\log n),
%$$
%then $s_n$ is summable by the $P$-method.
\end{theorem}
\begin{theorem}\label{t3} Let $P$ have all moments finite. If $s_n=O(n^{1/2}\log n)$, then $s_n\rightarrow s$ $(P)$ implies $s_n\rightarrow s$ $(\ell)$.
\end{theorem}
As Riesz means are central for us, it is as well to note here that it is not the function $\lambda$ in the Riesz mean $R(\lambda(n),1)$ that is relevant as such: what matters is the {\it order of magnitude of its logarithm}.  This is the content of the {\it second consistency theorem for Riesz means}; see e.g.~\cite{CM}, Ch. II and the references cited there.
\subsection{Proofs}
{\it Proof of Theorem \ref{e1}.} We adapt the approach of~\cite{Agn}. Recall first that for a transformation $\sigma_n=\sum_{k=0}^\infty a_{nk}s_k$, which assigns the value $\lim_{n\rightarrow\infty}\sigma_n$ to the sequence $\{s_n\}_{n=0}^\infty$, to be regular, by the Silverman-Toeplitz theorem it is necessary and sufficient for the following three conditions to hold: \\ (i) $\sum_{k=0}^\infty|a_{nk}|$ is bounded for all $n$, \\
(ii) $\lim_{n\rightarrow\infty}a_{nk}=0$, for each $k$, \\
(iii) $\lim_{n\rightarrow\infty}\sum_{k=1}^\infty a_{nk}=1$. Let $d_0 := s_0$, $d_1 := s_1$, and
\begin{eqnarray}
d_n:=\frac{1}{\log n}\sum_{n^{1/\lambda}<k\leq n}\frac{s_k}{k+1},\quad n\geq2.\nonumber
\end{eqnarray}
$\ell\Rightarrow\mathcal{L}(\lambda)$.  Let the sequence $\{t_n\}_{n=0}^\infty$ converge to $s$. It is clear that $d_0=t_0$, $d_1=t_1$ and
\begin{eqnarray}
d_n=t_n-(\log[n^{1/\lambda}])(\log n)^{-1}t_{[n^{1/\lambda}]},\quad n\geq 2.\nonumber
\end{eqnarray}
Thus $\mathcal{L}(\lambda)$ can be seen as a transformation of the sequence $\{t_n\}_{n=0}^\infty$. Moreover, it is a regular transformation. Indeed, conditions (i) and (ii) are clearly satisfied. For condition (iii) we need to show that
\begin{eqnarray}
\lim_{n\rightarrow\infty}\left[1-\frac{\log[n^{1/\lambda}]}{\log n}\right]=(1-\lambda^{-1}).\label{limit}
\end{eqnarray}
This follows by direct calculation, or (as the function $\log x$ is self-neglecting,~\cite{BinGT} \S 2.11,~\cite{BinO}) from the proof of Lemma in~\cite{BinG4}.\\
% that $\log[n^{1/\lambda}]\sim\log(n^{1/\lambda})=\lambda^{-1}\log n$, which proves (\ref{limit}) and thus $l\Rightarrow\mathcal{L}(\lambda)$.\\
$\mathcal{L}(\lambda)\Rightarrow \ell$. Let the sequence $\{d_n\}_{n=0}^\infty$ converge to $(1-\lambda^{-1})s$. Note that $t_0=d_0$, $t_1=d_1$, and
\begin{eqnarray}
t_n=d_n+\frac{\log[n^{1/\lambda}]}{\log n}d_{[n^{1/\lambda}]}+\frac{\log[[n^{1/\lambda}]^{1/\lambda}]}{\log n}d_{[[n^{1/\lambda}]^{1/\lambda}]}+....\label{tt}
\end{eqnarray}
Of course, the expansion (\ref{tt}) only contains a finite number of terms, and the final term depends on $n$ and $\lambda$. It is important to note that the coefficients in front of the elements $d_k$ of (\ref{tt}) have a certain pattern. The $\ell$-method can thus be seen as a transformation of the sequence $\{d_n\}_{n=0}^\infty$; we now show that it is regular. Condition (i) is clearly satisfied, and also (ii) since the coefficients in front of the elements $d_k$ are either zero or are decreasing with $n$. From (\ref{tt}) it is clear that condition (iii) follows from $\lim_{n\rightarrow\infty}\sum_{k=0}^n\lambda^{-k}=(1-\lambda^{-1})^{-1}$. \hfill$\Box$\\
%{\it Proof of Corollary \ref{e3}.} This follows from $R(\log n)\Leftrightarrow l\Leftrightarrow\mathcal{L}(\lambda)$, where $R(\log n)\Leftrightarrow l$ is Theorem 37 in~\cite{Har}. \hfill{$\Box$}\\

{\it Proof of Theorem \ref{u}.} Let (\ref{Ln}) hold for all $\lambda>1$. We can write (\ref{Ln}) as
\begin{eqnarray}
\frac{U(n)-U(n^{1/\lambda})}{\log n}\rightarrow(1-\lambda^{-1})s,\quad \forall \lambda>1,\label{U}
\end{eqnarray}
which clearly holds even for $\lambda=1$. By introducing $\alpha_n=\lambda^{-1}\log n$ and $V := U\circ \exp$, we can rewrite (\ref{U}) as
\begin{eqnarray}
\frac{V(\lambda \alpha_n)-V(\alpha_n)}{\alpha_n}\rightarrow(\lambda-1)s,\quad \forall \lambda\geq1.\label{V}
\end{eqnarray}
Since the linear function $x$ is regularly varying of index $1$, the function $V$ belongs to the de Haan class $\Pi_1$ (see Chapter 3 of~\cite{BinGT}). Hence the proof of the local uniformity follows from the proof of Theorem 3.1.16 of~\cite{BinGT} by using $\alpha_n$ instead of a continuous variable.\\
% Proof of (ii) cut-and-pasted to after end document

{\it Proof of Corollary \ref{e4}}. By the continuous version of Theorem \ref{u}, $\mathcal{L}_x$ for all $\lambda>1$ can be written as
\begin{eqnarray}
\frac{V(\lambda y)-V(y)}{y}\rightarrow(\lambda-1)s,\quad \forall \lambda\geq1,\label{V2}
\end{eqnarray}
where $y=\lambda^{-1}\log x$. From Theorem 3.2.7 of~\cite{BinGT} it now follows that (\ref{V2}) holds if and only if $s_n\rightarrow s$ $(\ell_x)$.\hfill{$\Box$}\\

{\it Proof of Theorem \ref{e6}.} This follows from (\ref{V2}) and Theorem 3.8.4 of~\cite{BinGT}.\hfill{$\Box$}\\

{\it Proof of Corollary \ref{e7}.} This follows from $\ell\Leftrightarrow\mathcal{L}_x(\lambda)$ for all $\lambda>1$, and taking $U_2 = 0$ in Theorem \ref{e6} (i). \hfill{$\Box$}\\

{\it Proof of Theorem \ref{t1}.} Recall the Hardy-Littlewood-Karamata theorem for the Laplace-Stieltjes transform (\cite{Bin6}, or \cite{BinGT}, \S 1.7): if the Laplace-Stieltjes transform of the function $G(x)$ is defined as $\hat{G}(s)=\int_o^\infty e^{-sx}dG(x)$, $H$ is slowly varying, and $\rho\geq 0$,
\begin{eqnarray}
G(x)\sim x^\rho H(x)/\Gamma(1+\rho),\quad (x\rightarrow\infty) \label{G}
\end{eqnarray}
is equivalent to
\begin{eqnarray}
\hat{G}(s)\sim H(1/s)/s^\rho,\quad (s\downarrow 0) \label{Gs}
\end{eqnarray}
if and only if
\begin{eqnarray}
\lim_{\lambda\downarrow 1}\liminf_{x\rightarrow \infty}\inf_{t\in[1,\lambda]}\frac{G(xt)-G(x)}{x^\rho H(x)}\geq 0.\label{HLK}
\end{eqnarray}
For $G(x)=U(x)$, $\rho=0$, $H(x)=\log x$, (\ref{G}) is in fact $l_x$, which is equivalent to $l$ (see Corollary \ref{e4}).  Then (\ref{Gs}) is equivalent to $L$.\hfill{$\Box$}\\

{\it Proof of Theorem \ref{t2}.} Part (i) follows from Theorem 1 in~\cite{Bin4}, and the fact that $C_1\Rightarrow \ell$, with $C_1$ being the Ces\`{a}ro summability method (see pp. 59 of~\cite{Har}). Part (ii) follows from
Theorem 3 of~\cite{Bin4}, and the fact that $C_{2r}\Rightarrow A\Rightarrow L$ by Theorem 55 of~\cite{Har} and Theorem 8 of~\cite{Ish1}, respectively (here $A$ denotes the Abel summability method).
For part (iii), note that with condition (\ref{P}) instead of $t_n= s+o(1/\log n)$ the proof of Theorem 2 in~\cite{Ish1} gives $(n+1)^{-1}\sum_{i=0}^n(s_i+\varepsilon_i)=s+o(n^{-1/2})$, with $\varepsilon_n\rightarrow 0$, and the conclusion follows from Proposition 2 of~\cite{Bin3}.\hfill{$\Box$}\\

{\it Proof of Theorem \ref{t3}.} The sequence $s_n=O(n^{1/2}\log n)$ is also $s_n=O(n^2)$ and thus of finite order. Corollary 2 in~\cite{Bin3} gives $P\Leftrightarrow V_1$, where $V_1$ is the Valiron method with parameter $1$ (see~\cite{Bin3} for the definition). Corollary 3 in~\cite{Bin3} gives $V_1\Leftrightarrow B$ with $B$ being the Borel method. Theorem in~\cite{Kwee3} gives $B\Rightarrow R(\log n)$.\hfill{$\Box$}\\

\subsection{Ordinary convergence}

The prototypical Tauberian theorem is one passing from convergence of an average (a logarithmic mean here) to ordinary convergence.  Necessary and sufficient Tauberian conditions for logarithmic convergence to imply ordinary convergence $s(.) \to c$ have recently been given by M\'oricz~\cite{Mor}.  These are of one-sided type, and involve logarithmic moving averages, as here:
$$
\limsup_{\l \downarrow 1} \liminf_{x \to \infty} \frac{1}{(\l - 1) \log x} \int_x^{x^{\l}} (s(u) - s(x)) du/u \geq 0,
$$
$$
\limsup_{\l \uparrow 1} \liminf_{x \to \infty} \frac{1}{(1 - \l) \log x} \int_{x^{\l}}^x (s(x) - s(u)) du/u \geq 0,
$$
(with one `$\liminf \limsup |...| = 0$' condition in the complex case).  These may be compared with the `$\liminf \liminf ... \geq 0$' necessary and sufficient (Tauberian) conditions for {\it Frullani integrals} in~\cite{BinG2} \S 6 (cf.~\cite{BinGT} \S 1.6.4).

\subsection{Remarks}
The analog of Theorem \ref{e1} for the Ces\`aro method is given in~\cite{Agn}. For a certain class of moving averages and their equivalence with Riesz means, see~\cite{BinG4}. It is interesting to note that Corollary \ref{e3} establishes an equivalence relation between a moving average method and a Riesz mean which is not included in the general result of~\cite{BinG4}.\\

The analog of Theorem \ref{t1} for the Ces\`aro and Abel methods is given in~\cite{Bin6}; for Euler, Borel, and $R(e^{\sqrt{n}})$ in~\cite{Bin6},~\cite{Bin2} (and Karamata-Stirling in~\cite{Bin7}); for the Valiron $V_\beta$ and Riesz $R(e^{n^{1-\beta}})$, $0<\beta<1$, in~\cite{BinT}.\\

There are a number of Tauberian and other theorems for the $\ell$ and $L$ methods (see, for example,~\cite{Ish0} - \cite{Ish3},~\cite{Ko1},~\cite{Ko2},~\cite{Kwee1},~\cite{Sit}).
\section{Law of large numbers and the Lambert W function}
\subsection{Results} Let $W(z)$ denote the Lambert W function, which is defined as the solution to the equation $z=W(z)e^{W(z)}$, $z\in\mathds{C}$ (see, e.g.~\cite{W}). Also let
\begin{eqnarray}
\phi(x):=(x+1)\log(x+1).\nonumber
\end{eqnarray}
\begin{lemma}\label{La} The function $W(z)$, $0\leq z\in\mathds{R}$, is subadditive.
\end{lemma}
\begin{theorem}\label{LLN} Let $X,X_1,X_2,...$, be a sequence of $i.i.d.$ random variables, and $m_k\equiv\mathds{E}[X_k\mathds{1}_{\{|X_k|\leq \phi(k+1)\}}]$. The following statements are equivalent:\\

(i) $\mathds{E}\left[e^{W(|X|)}\right] < \infty$, i.e. $\mathds{E}\left[\frac{|X|}{\log |X| \wedge 1}\right]<\infty$, or  $\mathds{E}\left[\frac{|X|}{1 + {\log}_+ |X|}\right]<\infty$, \\
% Simpler to state Th. 8 (ii) as

(ii) $X_n/(n\log n) \to 0$ $a.s.$ ($n \to \infty$),\\

(iii) $(X_n-m_n)\rightarrow 0$ $a.s.$ (l),\\

(iv) $(X_n-m_n)\rightarrow 0$ $a.s.$ (L),\\

(v) $(X_n-m_n)\rightarrow 0$ $a.s.$ (R($\log n$), or $\mathcal{L}(\lambda)$,  or $\mathcal{L}_x(\lambda)$,  or $\ell_x$),\\
%(vi) $(X_n-m_n)\rightarrow 0$ $a.s.$ ($\mathcal{L}(\lambda)$, \ \hbox{or) \ R(\log n), \ \hbox{or} \ ,\\
%(vii) $(X_n-m_n)\rightarrow 0$ $a.s.$ ($\mathcal{L}_x(\lambda)$),\\
%(viii) $(X_n-m_n)\rightarrow 0$ $a.s.$ ($l_x$).

(vi) $\frac{1}{\phi(n)}\sum_{1\leq i\leq n} (X_i-m_i)\rightarrow 0$ $a.s.$ ($n \to \infty$), or with $x$ instead of $n$,\\

(vii) $\frac{1} {\phi(n)}\sum_{\phi^{\leftarrow}(\beta^{-1}\phi(n))<i\leq n}(X_i-m_i)\rightarrow 0$ $a.s.$ ($n \to \infty$), or with $x$ instead of $n$, $\forall \beta>1$.\\

(viii) $\sum_1^\infty n^{-1}\mathds{P}[|\sum_{1\leq i\leq n}(X_i-m_{i+n/(\gamma-1)})|>\phi(n/(\gamma-1))\epsilon]<\infty$ $\forall\epsilon>0$ and $\forall \gamma>1$,\\

(ix) $\sum_1^\infty n^{-1}\mathds{P}[\max_{1\leq k\leq n}|\sum_{1\leq i\leq k}(X_i-m_{i+n/(\gamma-1)})|>\phi(n/(\gamma-1))\epsilon]<\infty$ $\forall\epsilon>0$ and $\forall \gamma>1$.

%If $\frac{1}{\phi(n)}\sum_{1\leq i\leq n} m_i\rightarrow 0$ ($n \to \infty$), and for $\gamma=2$, the statements $(viii)$ and $(ix)$ take the following simpler form:\\

%$(viii^*)$ $\sum_{n=1}^\infty n^{-1}\mathds{P}[|\sum_{i=1}^nX_i|>\phi(n)\epsilon]<\infty$   $\forall \epsilon >0$,\\

%$(ix^*)$ $\sum_{n=1}^\infty n^{-1}\mathds{P}[\max_{1\leq k\leq n}|\sum_{i=1}^kX_i|>\phi(n)\epsilon]<\infty$   $\forall \epsilon >0$.
\end{theorem}
\subsection{Proofs}
{\it Proof of Lemma \ref{La}}. Let $0\leq a\in\mathds{R}$ and define $F(z):=W(z+a)-W(z)$. Then
\begin{eqnarray}
\frac{dF(z)}{dz}&=&\frac{e^{-W(z+a)}}{1+W(z+a)}-\frac{e^{-W(z)}}{1+W(z)}\nonumber\\
&=&\frac{e^{-W(z+a)}[1+W(z)]-e^{-W(z)}[1+W(z+a)]}{[1+W(z+a)][1+W(z)]}<0,\nonumber
\end{eqnarray}
since $W(z)$ is an increasing function for $z\geq0$. This means that $F(z)$ is a decreasing function with its maximum value of $F(0)=W(a)$, which implies $W(z+a)-W(z)\leq W(a)$.\hfill{$\Box$}\\

{\it Proof of Theorem \ref{LLN}.} The asymptotics of the equation $z = W e^W$ for the Lambert W function are considered in detail by de Bruijn (\cite{Bru}, \S 2.4). From this,
$$
e^W \sim z/\log z \qquad (z \to \infty).
$$
So as $\mathds{E}[e^{W(|X|)}] < \infty$ is restrictive only for large values of $e^{W(|X|)}$, the three moment conditions in (i) are equivalent.\\
It is clear from the previous section that $(iii)\Leftrightarrow(v)$. Note that the $\ell$-method is equivalent to
\begin{eqnarray}
\frac{1}{\log (n+1)}\sum_{i=1}^n\frac{s_i}{i+1}\rightarrow s, \quad (n\rightarrow\infty).\nonumber
\end{eqnarray}
The inverse of $\phi(x):= (x+1)\log(x+1)$ is $\phi^{-1}(x)=-1+e^{W(x)}$. The equivalence $(i)\Leftrightarrow(iii)$ now follows from a result of Jajte ~\cite{jajte}, from the proof of which $(iii) \Rightarrow (ii) \Rightarrow (i)$, so $(i) - (iii)$ are equivalent.\\
The equivalence of $(i)$ and $(iv)$ is contained in a result of Kiesel (\cite{Kie}, Th., p.196; also p.198, 203).  As Kiesel points out, the prototypical results here are those of Chow~\cite{Chow} and Lai~\cite{Lai}; cf.~\cite{CL}.  The proof is via symmetrization inequalities and a Borel-Cantelli argument. The following is a slightly different proof of $(i)\Leftrightarrow(iv)$ as compared to~\cite{Kie}. Due to the Abelian fact that always $\ell\Rightarrow L$, we have $(i)\Rightarrow(iii)\Rightarrow(iv)$. $(iv)$ implies
\begin{eqnarray}
\frac{1}{\log m}\sum_{k=1}^\infty \frac{X_k^s}{k+1}e^{-(k+1)m^{-1}}=0,\quad a.s.,\nonumber
\end{eqnarray}
where $X_k^s=X_k-X_k'$, and $\{X_n\}_{n=1}^\infty$ and $\{X'_n\}_{n=1}^\infty$ are i.i.d.. We define
\begin{eqnarray}
\widetilde{X}_m:=\frac{1}{\log m}\sum_{k=1}^m \frac{X_k^s}{k+1}e^{-(k+1)m^{-1}}, \quad \widehat{X}_m:=\frac{1}{\log m}\sum_{k=m+1}^\infty \frac{X_k^s}{k+1}e^{-(k+1)m^{-1}}.\nonumber
\end{eqnarray}
Then $\widetilde{X}_k$ and $\widehat{X}_m$ are independent and symmetric. Since $\widetilde{X}_m+\widehat{X}_m\rightarrow 0$, a.s., we also have convergence in probability. From the L\'{e}vy inequality (Lemma 2 in V.5 of~\cite{Fel}),
$\widehat{X}_m\rightarrow 0$ in probability. Since $(\widetilde{X}_1,...,\widetilde{X}_m)$ and $\widehat{X}_m$ are also independent, Lemma 3 of~\cite{CL} gives $\widetilde{X}_m\rightarrow 0$, $a.s.$. Writing
\begin{eqnarray}
\widetilde{X}_m=\frac{1}{\log m}\sum_{k=1}^{m-1}\frac{X_k^s}{k+1}e^{-(k+1)m^{-1}}+\frac{1}{\log m}\frac{X_m^s}{(m+1)}e^{-(m+1)m^{-1}}\nonumber
\end{eqnarray}
and repeating the previous argument, we obtain $[(m+1)\log(m+1)]^{-1}X_m^s\rightarrow 0$, $a.s.$. By the Borel-Cantelli lemma, and the weak symmetrisation inequalities (pp. 257 of~\cite{Loe}),
\begin{eqnarray}
\frac{1}{2}\sum_{k=1}^\infty\mathds{P}\left[-1+e^{W(|X-\mu_x|)}\geq k\right]=\frac{1}{2}\sum_{k=1}^\infty\mathds{P}\left[|X-\mu_x|\geq (k+1)\log(k+1)\right]\nonumber\\
\leq\sum_{k=1}^\infty\mathds{P}\left[|X^s|\geq (k+1)\log(k+1)\right]<\infty,\nonumber
\end{eqnarray}
with $\mu_x$ the median of $X$, and $X^s=X-X'$, with $X$ and $X'$ i.i.d. This and $W$ subadditive give
\begin{eqnarray}
\mathds{E}\left[e^{W(|X|)}\right]\leq\mathds{E}\left[e^{W(|X-\mu_x|+|\mu_x|)}\right]\leq e^{W(|\mu_X|)}\mathds{E}\left[e^{W(|X-\mu_x|)}\right]<\infty.\nonumber
\end{eqnarray}
The equivalence $(i)\Leftrightarrow (vi)$ when the convergence in $(vi)$ is through the integers $n$ follows from Jajte's result~\cite{jajte}. The equivalence $(vi)\Leftrightarrow (vii)$ in this case can be obtained in a similar way to the equivalence $\ell\Leftrightarrow\mathcal{L}(\lambda)$ of the previous section. The equivalence relations for the continuous  variable $x$ are obtained in a same way as in Corollary \ref{e4}.\\
We now show $(vi)\Leftrightarrow (viii)$. Since $\phi(x)$ is a regularly varying function of index 1, it follows from Theorem 3.2.7 of~\cite{BinGT} that $(vi)$, when the convergence is through the continuous variable $x$, is equivalent with
\begin{eqnarray}
\frac{1}{\phi(x)}\sum_{x<i\leq \gamma x}(X_i-m_i)\rightarrow 0\quad a.s.\quad (x\rightarrow\infty)\quad \forall\gamma>1,\label{BK}
\end{eqnarray}
and the convergence is locally uniform on $(1,\infty)$. By introducing the variable $y:=\log x/\log \gamma$, i.e. $x=\gamma^y$, we rewrite (\ref{BK}) as
\begin{eqnarray}
\frac{1}{\phi(\gamma^y)}\sum_{\gamma^y<i\leq \gamma^{y+1}}(X_i-m_i)\rightarrow 0\quad a.s.\quad (y\rightarrow\infty)\quad \forall\gamma>1,\nonumber
\end{eqnarray}
which is equivalent to
\begin{eqnarray}
Z_n:=\frac{1}{\phi(\gamma^n)}\sum_{\gamma^n<i\leq \gamma^{n+1}}(X_i-m_i)\rightarrow 0\quad a.s.\quad (n\rightarrow\infty)\quad \forall\gamma>1.\label{z}
\end{eqnarray}
By the independence of random variables $X_i$ and the non-overlapping of the defining summations for different $n$, the $Z_n$ are independent. Thus by the Borel-Cantelli lemma, (\ref{z}) is equivalent to
\begin{eqnarray}
\sum_{i=1}^\infty \mathds{P}[|Z_n|>\epsilon]<\infty\quad\forall\epsilon>0.\label{BC}
\end{eqnarray}
Since the law of $\sum_{\gamma^n<i\leq \gamma^{n+1}}X_i$ is the same to that of $\sum_{1\leq i\leq \gamma^{n+1}-\gamma^n}X_i$, we can write (\ref{BC}) as
\begin{eqnarray}
&&\sum_{i=1}^\infty \mathds{P}\left[\left|\frac{1}{\phi(\gamma^n)}\sum_{\gamma^n<i\leq \gamma^{n+1}}(X_i-m_i)\right|>\epsilon\right]\nonumber\\
=&&\sum_{i=1}^\infty \mathds{P}\left[\left|\sum_{1\leq i\leq \gamma^{n}(\gamma-1)}(X_i-m_{\gamma^n+i})\right|>\phi(\gamma^n)\epsilon\right]<\infty\quad\forall\epsilon>0,\nonumber
\end{eqnarray}
which is equivalent to
\begin{eqnarray}
\int^\infty \mathds{P}\left[\left|\sum_{1\leq i\leq \gamma^{t}(\gamma-1)}(X_i-m_{\gamma^t+i})\right|>\phi(\gamma^t)\epsilon\right]dt<\infty\quad\forall\epsilon>0.\nonumber
\end{eqnarray}
After the change of variable $\gamma^{t}(\gamma-1):=u$ we obtain the equivalent inequality
\begin{eqnarray}
\int_1^\infty \mathds{P}\left[\left|\sum_{1\leq i\leq u}(X_i-m_{i+u/(\gamma-1)})\right|>\phi(u/(\gamma-1))\epsilon\right]\frac{1}{u\log\gamma}du<\infty\quad\forall\epsilon>0,\nonumber
\end{eqnarray}
which is equivalent to $(viii)$. The proof of $(vi)\Leftrightarrow(ix)$ follows in the same way, but instead of (\ref{BK}) the starting point is its equivalent (due to local uniformity)
\begin{eqnarray}
\frac{1}{\phi(x)}\max_{x<k\leq \gamma x}\left|\sum_{x<i\leq k}(X_i-m_i)\right|\rightarrow 0\quad a.s.\quad (x\rightarrow\infty)\quad \forall\gamma>1.\nonumber
\end{eqnarray}
\hfill{$\Box$}
\subsection{Remarks}
We call the moment condition in Theorem \ref{LLN}(i) the {\it $L/\log L$ condition}.  It expresses a little {\it less} integrability than in Kolmogorov's strong law of large
numbers (SLLN) of 1933, where $X \in L_1$, i.e. $E[|X|] < \infty$, is necessary and sufficient, and correspondingly one needs to divide by $n \log n$ rather than $n$ in Theorem \ref{LLN} (ii).  By contrast, the {\it $L \log L$ condition} $\mathds{E}[|X| {\log}_+ |X|] < \infty$, where $X$ has a little {\it more} integrability than in SLLN, occurs in various connections, e.g. branching processes and martingale theory.  See e.g.~\cite{AthN} I.C.10,~\cite{Nev} Prop. IV.2.10, 11. \\
\i With $\mu := \mathds{E}[X]$ where this exists ($X \in L_1$, $\mathds{E}[|X|] < \infty$), the crucial role of SLLN (see \S 4.1 below) is revealed in a {\it discontinuity in the functional form} regarding the Ces\`aro family $C_{\alpha}$ of summability methods ($\alpha > 0$: see~\cite{Har}, V, VI).  For $\alpha \geq 1$, the $C_{\alpha}$ are equivalent here:
$$
X \in L_1 \quad \Leftrightarrow \quad X_n \to \mu \quad a.s. \quad (C_1) \quad \Leftrightarrow X_n \to \mu \quad a.s. \quad (C_{\alpha}) \quad (\alpha \geq 1).
$$
But for $\alpha \leq 1$, the integrability required here is the more stringent $L_{1/\alpha}$:
$$
X \in L_{1/\alpha} \quad \Leftrightarrow \quad X_n \to \mu \quad a.s. \quad (R_{1/\alpha}) \quad \Leftrightarrow X_n \to \mu \quad a.s. \quad (C_{\alpha}) \quad (0 < \alpha \leq 1),
$$
with $R_p$ ($p := 1/\alpha$) the {\it Riesz mean} studied in~\cite{Bin5}.  The case $\alpha = \frac{1}{2}$, $p = 2$ ($X \in L_2$, finite variance) is particularly important, as it occurs in connection with the Euler and Borel summability methods, as in the papers of Chow~\cite{Chow} and Lai~\cite{Lai}. One can extend to more general moment conditions $\mathds{E}[\psi(|X|)] < \infty$, as in~\cite{BinG4}. \\

Our proof of Theorem \ref{LLN} above is different in detail from those of these sources, being based on the {\it Kiesel-Jajte} centring method of~\cite{Kie} and~\cite{jajte}.  But it shares a number of aspects with the Chow-Lai arguments above.  Our aim in setting out the proof has been to keep the probability and summability (or analysis) aspects separate. \\

While in the results above it matters crucially whether one has less or more integrability than in SLLN, two classical results (of a different nature to those above) enable one to handle both together.  These are the Marcinkiewicz-Zygmund SLLN~\cite{MarZ} and the Baum-Katz LLN~\cite{BauK}; see e.g.~\cite{Bin8} for background and references.\\

We use $\phi(n) = (n+1)\log (n+1)$ rather than $\phi(n) =n \log n$ as in~\cite{jajte} to avoid a slip made there:
%, a slightly different version of the $l$-method was considered with  instead of our  . The following claim is made there
\begin{eqnarray}
\mathds{E}[|X|^\alpha]\leq\mathds{E}[\phi^{-1}(|X|)]\leq\mathds{E}[|X|],\quad 0<\alpha<1,\nonumber
\end{eqnarray}
is claimed (but for $0\leq X\leq 1/2$ the right inequality does not hold since $\mathds{E}[\phi^{-1}(|X|)]\geq 1$ and $\mathds{E}[|X|]\leq 1/2$).\\

The summability methods of Theorem \ref{LLN} $(vi)$ and $(vii)$ are {\it non-regular}, giving an equivalence between regular and non-regular summability methods in this probabilistic setting. Parts $(viii)$ and $(ix)$ are the LLN of Baum-Katz type for the logarithmic summability methods. In~\cite{Lai2} Lai gives another result of Baum-Katz type in this context, whereas~\cite{BinG4} has a Baum-Katz type LLN for a general class of moving averages.\\

There are many analogs to parts $(i)- (iv)$ of Theorem \ref{LLN} for other summability methods. Few of them are: for the Ces\`{a}ro and Abel methods~\cite{Lai}; for the Borel and Euler methods~\cite{Chow}; for Karamata-Stirling and Jakimovski methods~\cite{BinS}, for the $P$ and Valiron methods~\cite{BinM}.\\

Kiesel~\cite{Kie2} considers analogues of the law of the iterated logarithm (LIL) for power-series methods with weights $p_n$ regularly
varying with index $\alpha > - 1/2$ (or $\alpha > -1$ in the appendix). His interesting results do not apply to the logarithmic case here with
$\alpha = -1$.\\

Some further problems that could be considered in this direction are: $\phi$-mixing case as in~\cite{Bin8}, Banach space case as in~\cite{Bin8}, negatively associated random variables as in~\cite{BinN}, and the `law of the single logarithm' (LSL) as in~\cite{Bin5},~\cite{GutJS}. \\
\section{Applications}
We restrict to three, one in probability theory, two in analytic number theory.
\subsection{The almost-sure central limit theorem (ASCLT)}.
Recall the two central pillars of probability theory: with $X, X_1, \ldots, X_n$ i.i.d. random variables with mean $\mu$, $S_n := \sum_1^n X_k$,
$$
S_n/n \to \mu = E[X] \qq (n \to \infty) \qq a.s.
$$
(SLLN), and if the $X_n$ have variance ${\s}^2$,
$$
P((S_n - n \mu)/(\s \sqrt{n}) \leq x) \to \Phi(x) := \int_{-\infty}^x \frac{e^{- \half u^2}}{\sqrt{2 \pi}} du \q (n \to \infty) \q \forall \ x \in \R
$$
 (central limit theorem, CLT).  The ASCLT combines aspects of both: taking $\mu = 0$, $\s = 1$ for simplicity,
 $$
 \frac{1}{\log n} \sum_{k=1}^n I(S_k/\sqrt{k} \leq x)/k \to \Phi(x) \q (n \to \infty) \q a.s. \q \forall \ x \in \R.
 $$
 \i Results of this type can be traced back to L\'evy in 1937 (\cite{Lev}, p.270); see~\cite{Bin96} for the difference between logarithmic means (a.s. convergence) and Ces\`aro means (convergence in distribution to the arc-sine law). It was proved in the form above by Brosamler and Schatte independently in 1988.  There have been many extensions.  We refer for details and references to the surveys by Peligrad and R\'ev\'esz~\cite{PelR} in 1991, Berkes~\cite{Ber} in 1998, Berkes et al.~\cite{BerCCM} in 2002.  Versions of the ASCLT with non-standard weights are considered (that is, the extent to which the logarithmic method may be generalised here).  It turns out that, if limits as $n \to \infty$ are allowed to omit an exceptional set of logarithmic density 0, much of probability limit theory changes: for instance, the weak (convergence in probability) and strong laws of large numbers become equivalent. \\
 \i Such exceptional sets of logarithmic density 0 occur in regular variation in analysis; see e.g.~\cite{BinGT}, \S 2.9. \\

 \subsection{Number-theoretic densities}
 \i A subset $A \subset \N$ has lower and upper {\it arithmetic densities}
 $$
 \underline{\bf d} A := \liminf_{n \to \infty} \frac{1}{n} | A \cap \{1,2, \ldots, n \}|, \q
 \overline{\bf d} A := \limsup_{n \to \infty} \frac{1}{n} | A \cap \{1,2, \ldots, n \}|,
 $$
 and {\it arithmetic density} ${\bf d} A$ their common value when these are equal.  The upper and lower {\it logarithmic densities} are defined by
 $$
 \underline{\delta} A := \liminf_{n \to \infty} \frac{1}{\log n} \sum_{k \leq n, k \in A} 1/k, \q
 \overline{\delta} A := \limsup_{n \to \infty} \frac{1}{\log n}  \sum_{k \leq n, k \in A} 1/k,
 $$
with the {\it logarithmic density} $\delta A$ their common value when equal.  Then
$$
\underline{\bf d} A \leq \underline{\delta} A \leq \overline{\delta} A \leq \overline{\bf d} A:
$$
existence of ${\bf d}$ implies that of $\delta$, but not conversely, so $\delta$ extends ${\bf d}$.  See e.g. Tenenbaum (\cite{Ten}, III.1.2). \\
\i For $\s > 1$ and $\z(\s) := \sum_1^{\infty} 1/n^{\s}$ the Riemann zeta function, write
$$
P_{\s}(A) := {\z(\s)}^{-1} \sum_{n \in A} 1/n^{\s}.
$$
If this has a limit as $\s \downarrow 1$, this is called the {\it analytic} (or {\it Dirichlet}) {\it density} of $A$ (equivalently, as $\z$ has a simple pole at 1 of residue 1, one can use $(\s - 1) \sum_{n \in A} 1/n^{\s}$ here).  Then (\cite{Ten}, III.1.2) $A$ has analytic density iff it has logarithmic density, and the two are equal. \\
\i One classical application is to Dirichlet's theorem of 1837 on primes in an arithmetic progression $a + nb$ (for $a$, $b$ coprime).  This can be shown to have analytic density $1/\phi(b)$ (for $\phi$ the Euler totient function), proving Dirichlet's theorem that there are infinitely many such primes.  Indeed, they have arithmetic density $1/\phi(b)$ (so the primes are equally distributed between the residue classes), but this is harder (de la Vall\'ee Poussin in 1896; see e.g. Burris (\cite{Bur}, 9.2)). \\
\i Further applications of logarithmic density in number theory are given in Halberstam and Roth (\cite{HalR}, V.5).

\subsection{The Prime Number Theorem (PNT)}

With $\Lambda$ the von Mangoldt function, $\Lambda(n) := \log p$ if $n = p^m$ is a prime power, 0 otherwise, one has
$$
\Lambda \to 1 \qquad ({\ell}_x): \qquad \sum_{n \leq x} \Lambda(n)/n \sim \log x \qquad (x \to \infty)                          \eqno(4.1)
$$
(and similarly with $\cal{L}(\lambda)$ for ${\ell}_x$, indeed uniformly on compact $\lambda$-sets in $(1,\infty)$, by Corollary 2 and Theorem 2). Though relevant to PNT, $\pi(x) \sim li(x) \sim x/\log x$ (with $\pi(x) := \sum_{p \leq x} 1$ the prime-counting function, $li(x) := \int_2^x dt/\log t$), this is far weaker.  From the stronger result
$$
\sum_{n \leq x} \Lambda(n)/n = \log x + O(1)                                                                                    \eqno(4.2)
$$
(see e.g.~\cite{HarW}, Theorem 424) one can deduce Tchebychev's result
$$
\liminf \pi(x)/li(x) \leq 1 \leq \limsup \pi(x)/li(x),
$$
so that if the limit exists, it is 1 (see e.g.~\cite{HarW}, Theorem 426,~\cite{Ten}, I.1.7).  But this is still far weaker than PNT, which is equivalent by elementary means to
$$
\sum_{n \leq x} \Lambda(n)/n = \log x + constant + o(1)                                                                         \eqno(4.3)
$$
(the constant is Euler's constant $\gamma$, but this is not needed here:~\cite{Ten}, I.3.6).  There is a similar hierarchy of results concerning the M\"obius function $\mu$.  For background and details, see~\cite{Dia}.

As $\sum_{n \leq x} \Lambda(n)/n = \sum_p \log p/p + \sum_{m \geq 2, p \leq x} \log p/p^m$ and the second sum is convergent as $x \to \infty$ (as one can easily check), $(4.3)$ is equivalent to
$$
\sum_{p \leq x} \log p/p = \log x + constant + o(1),
$$
another form of PNT (\cite{Dia}, (4.3)). The weaker $(4.2)$ is equivalent to a result of Mertens,
$$
\sum_{p \leq x} \log p/p = \log x + O(1)
$$
(\cite{Ten}, I.1.4, ~\cite{HarW}, Theorem 425).  See~\cite{BinI} for background here.

With $p_n$ the $n$th prime, a further equivalent form of PNT is $p_n \sim n \log n$ (\cite{Dia}, (4.6),~\cite{Ten}, I Ex. 5a), which is to PNT in the form $\pi(x) \sim x/\log x$ as (ii) is to (i) in Theorem \ref{LLN}.  Note in this regard that the equivalent moment conditions in Theorem \ref{LLN}(i) may be augmented by $E[\pi(|X|)] < \infty$, $E[li(|X|)] < \infty$. \\

%With $p_n$ the $n$th prime, a further equivalent form of PNT is $p_n \sim n \log n$ (\cite{Dia}, (4.6),~\cite{Ten}, I Ex. 5a). Since $\Lambda(n)=o(n\log n)$, by Corollary of~\cite{Sit}, we have $\Lambda(n)\rightarrow (\ell)$ if and only if $\Lambda(n)\rightarrow (L)$. We thus have the following PNT analog to Theorem \ref{LLN}.
%\begin{theorem} The following statements are true and equivalent:\\

%(i) $e^{W(x)}\sim \pi(x)\sim x/\log x$,\\

%(ii) $p_n \sim n \log n$,\\

%(iii) $\Lambda(n) \rightarrow 1$ $(\ell_x)$,\\

%(iv) $\Lambda(n) \rightarrow 1$ $(L)$,\\

%(v) $\Lambda(n) \rightarrow 1$ (R($\log n$), or $\mathcal{L}(\lambda)$,  or $\mathcal{L}_x(\lambda)$,  or $\ell$)
%\end{theorem}

\section*{Acknowledgement}
We thank the referee for his careful and informative report.

\end{document}